# EFFECTS OF ELECTRIC VEHICLE ADOPTION FOR STATE-WIDE INTERCITY TRIPS ON EMISSION SAVING AND ENERGY CONSUMPTION

Mohammadreza Kavianipour[a], Hamid Mozafari[a], Mehrnaz Ghamami[a*], Ali Zockaie[a], Robert Jackson[b]

[a] Michigan State University, 428 S Shaw Ln, East Lansing, MI 48824, USA
[b] Michigan Department of Environment, Great Lakes, and Energy. 525 West Allegan Street, Lansing, MI 48909-7973, US

[*]Corresponding author: ghamamim@msu.edu, Tel: +1 517 355 1288
Email addresses: kavianip@egr.msu.edu (Mohammadreza Kavianipour), mozafar1@msu.edu (Hamid Mozafari),  zockaiea@egr.msu.edu (Ali Zockaie), JacksonR20@michigan.gov (Robert Jackson)

## 1. ABSTRACT

Electric vehicles (EVs) mitigate fossil fuel dependency and reduce traffic emissions; hence they are considered sustainable alternatives to conventional vehicles. To support EVs intercity trips, a recent study proposed a charging infrastructure planning framework, which is implemented for the intercity network of Michigan considering different projected EV market shares of 2030. This study aims to estimate the emission reduction associated with the projected electrification rate and the proposed infrastructure for light-duty vehicles in the earlier study. To this end, a state-of-the-art emission estimation framework is proposed for state-wide intercity travels. The main contributions of the proposed framework include: 1) Incorporating a micro-emission estimation model for the simulated vehicle trajectories of the intercity network of Michigan; 2) Adjusting the micro-emission model results considering impacts of monthly travel demand and temperature variations, and heterogeneity of vehicles based on their make, model, and age; 3) Contrasting the required investment on charging infrastructure to support EVs intercity trips and the associated emission savings for a given market share. To this end, the proposed emission estimation framework is benchmarked against a traditional method based on vehicle miles traveled (VMT). Then, five different scenarios of EV adoption are explored to assess potential emission savings. The results suggest that the annual $CO_2$ savings due to the electrification of intercity travels range from 0.56 million-ton to 0.93 million-ton in 2030 for a projected 6 percent market share, which results in \$28M to \$46M $CO_2$ societal cost savings. These savings justify the investment in network electrification. Note that only 3.8 to 8.8 percent of the EV energy demand must be satisfied by the DC fast charger network proposed by the charging infrastructure planning framework. This requires 22.45 to 51.60 BWh annual energy consumptions associated with the estimated EV market share in Michigan for 2030.

Kavianipour et al.    3

## 2. INTRODUCTION

Transportation is a major source of emission. In 2018, 28% of total Greenhouse Gas (GHG) emission was originated from the transportation sector in the U.S. [1]. Air pollution has caused significant health problems in recent years. Many studies have associated air pollution with allergies, heart attacks, asthma, and lung cancer [2–5]. Long-term exposure to air pollution increases the mortality rate [6]. Poor air quality is responsible for five to ten percent of total premature death in the U.S. [7]. Therefore, various approaches have been implemented to reduce GHG emission including promoting non-motorized modes of transportation [8], eco-driving [9], congestion mitigation strategies such as pricing [8], using more fuel-efficient vehicles [10], and alternative fuel vehicles such as liquefied natural gas [11] or electric vehicles (EV) [12–15].

Macroscopic and microscopic models are implemented extensively to estimate vehicle emissions at different scales. The macroscopic models have been the dominant approach for estimating vehicle emissions for many years. They typically calculate emission per unit of distance based on the average travel speed [16]. While the average speed affects the emission generation rate significantly, the instantaneous fluctuations of speed are also a major contributor, which is ignored in macroscopic models. Thus, applying macroscopic models might be appropriate for large-scale applications [16]; however, it entails an approximation, which is explored in our study. IVE, MOBILE, and COPERT are the most commonly used models in this category [17].

On the other hand, the micro-emission modeling can be applied at both vehicular and network-level [18,19]. These models can capture the correlation between emission and vehicle dynamics, e.g., speed and acceleration. Several studies in the literature have collected vehicle emission data in real-world traffic conditions using instrumented vehicles to calibrate these emission models [16,20]. These models can estimate emissions for a road segment, or even capture the impacts of traffic-management policies and different driving patterns by considering the impacts of congestion and instantaneous speed fluctuations in estimating vehicle emission [21]. Comprehensive Modal Emissions Model (CMEM) [22] and International Vehicle Emissions (IVE) model [23] are two instances of micro-emission models.

Micro-emission models are usually coupled with traffic microsimulation models, which provide the required speed and acceleration profiles. These models simulate real-time individual driving behavior and trajectories. VISSIM [24], AIMSUN [25], and PARAMICS [26] are common traffic microsimulation tools, which are used to estimate vehicular emissions and fuel consumption [27]. As the micro-emission estimation approach is resource-intensive, it is not applicable in large-scale networks. Therefore, mesoscopic traffic simulation models are adopted for large-scale networks [28]. These models reduce the required computation time significantly compared to microscopic models, while improving the accuracy relative to macroscopic models. Several studies have used the mesoscopic approach to estimate GHG emission [20,28–31]. In this study, a micro-emission model coupled with a mesoscopic traffic simulator is compared with a macroscopic model as a benchmark.

The adverse weather conditions in Michigan during the winter seasons affect the battery performance, which yields in significant vehicle range reduction. The impact of operating temperature on vehicle emission is not widely studied as the majority of vehicle emission tests are carried out at normal temperatures, i.e., 70 ˚F [32]. The United States Environmental Protection Agency (EPA) suggests that in the states with harsh winter, the impact of cold temperature should be considered to prevent under-estimating the vehicle emissions [33]. Aligned with this recommendation, a recent study conducted experimental tests on vehicles at different ambient temperatures, namely 72 °F, 20 °F, and −1 °F, to demonstrate that the cold ambient temperature increases emission and fuel consumption [34]. This increase is further amplified as the vehicle ages [35,36], and has different amplitude in gasoline and diesel vehicles [37,38]. In a recent study, five gasoline vehicles and five diesel vehicles were tested at 20˚F and 73˚F using the World-Harmonized Light-duty Test Cycle (WLTC) [39]. When the temperature drops from 73˚F to 20˚F, NOx, HC, and CO2 productions escalate on average by 3.4, 1.5, and 1.15 times for diesel vehicles, while these proportions are 1.7, 6.5, and 1.11 for gasoline vehicles, respectively.

Vehicle emission also varies based on the vehicle make and model. The Environmental Protection Agency (EPA) presents emission rates for numerous vehicle models from 1984 to 2020 [40]. The vehicle emission rate varies by the vehicle specifications, such as make, model, manufactured year, and engine displacements. For instance, Barth et al. tested 24 types of gasoline vehicles, over the model years of 2001 to 2004, and showed significant variations over their average emission rate [41]. They used dynamometer testing and instrumented vehicles to measure the emissions, and vehicle mass, engine size, and transmission type were identified as the key factors affecting emission rate.

As mentioned earlier, EVs are proposed to reduce emissions in transportation networks due to their no-tailpipe emission. However, supporting intercity trips for these vehicles requires a significant capital investment on fast-charging infrastructure due to their limited range. Thus, the resulting emission savings from vehicle electrification needs to be compared with the required capital investment to assess this plan economically. We have explored the required charging infrastructure investment via an optimization framework in recent studies [14,15,42]. In these studies, the annual travel demand variations and seasonal weather conditions are considered state-of-the-art in developing charging infrastructure for EVs. Both of these factors are key factors in emission estimation as well. This study investigates the impact of intercity road network electrification on vehicle emission in Michigan and compares it with the required infrastructure investment and electricity consumption based on the earlier developed framework [14].

The key contributions of this study are as follows:
- Proposing a large-scale emission estimation framework coupling a micro-emission model and a mesoscopic traffic simulation tool, considering monthly travel demand variations, seasonal weather conditions, and vehicle type composition.
- Conducting a scenario analysis for vehicles replaced by EVs and resulting emission savings from the electrification of intercity trips in Michigan.
- Contrasting the emission savings with the required fast-charging infrastructure investment to support intercity trips for different EV market shares.

The pollutant considered in this study is $CO_2$, as the main running emission pollutant for light-duty vehicles. A mesoscopic dynamic simulation tool, DYNASMART [43], is used to simulate vehicle trajectories and provide the speed and acceleration profiles. Monthly demand factors are estimated to capture the seasonal variation of travel demand [15]. Then, an emission estimation model, proposed by Panis et al. [16], is incorporated to calculate emission from each trajectory. The estimated values are then adjusted based on different vehicle types traveling in Michigan, and the ambient temperature variations to estimate the vehicle emission at large-scale networks accurately. These results are compared with a macroscopic approach as a benchmark to demonstrate the importance of the developed framework in this study for intercity network emission estimation. Finally, it estimates the societal cost of vehicle emission reduction due to electrification and explores the trade-off between investing in the charging infrastructure and the savings of the societal costs considering the electric energy consumption of EVs.

## 3. METHODOLOGY

This section first presents the proposed research framework (Figure 1) by illustrating the relation between different research approach components, including simulation, monthly factors, and emission estimation. Then, it describes each component as follows. The first two components provide vehicle trajectories and Monthly origin-destination (OD) demand factors to the emission estimation component. The simulation component, highlighted in yellow, provides the zone to zone vehicle trajectories of the state-wide Michigan network, clustered to OD trajectories of the intercity network according to zones associated with each major city with a population of over 50,000. The monthly factor component incorporates the simulated link volumes for the base-month scenario and the real-world volumes captured during different months of year via continuous counting stations installed on Michigan highways to estimate demand factors.

Then, based on a micro-emission model, highlighted in blue, the emission associated with each OD pair in the intercity network, is estimated and modified for each month based on the monthly factors. Then another sub-module assigns a vehicle type (model, make, and year) randomly to each trajectory, which is highlighted in gray. This sub-module ensures that the distribution of the vehicle types matches the sample market share provided by the National Household Travel Survey (NHTS) [44]. The vehicle type emission adjustment factors, derived from EPA dataset [40], are then applied to adjust the trajectory-based emission. Finally, in the last sub-module, the weather adjustment factors are applied to the monthly emissions to estimate the final results. The remaining portion of this section describes the above modules and their associated sub-modules.

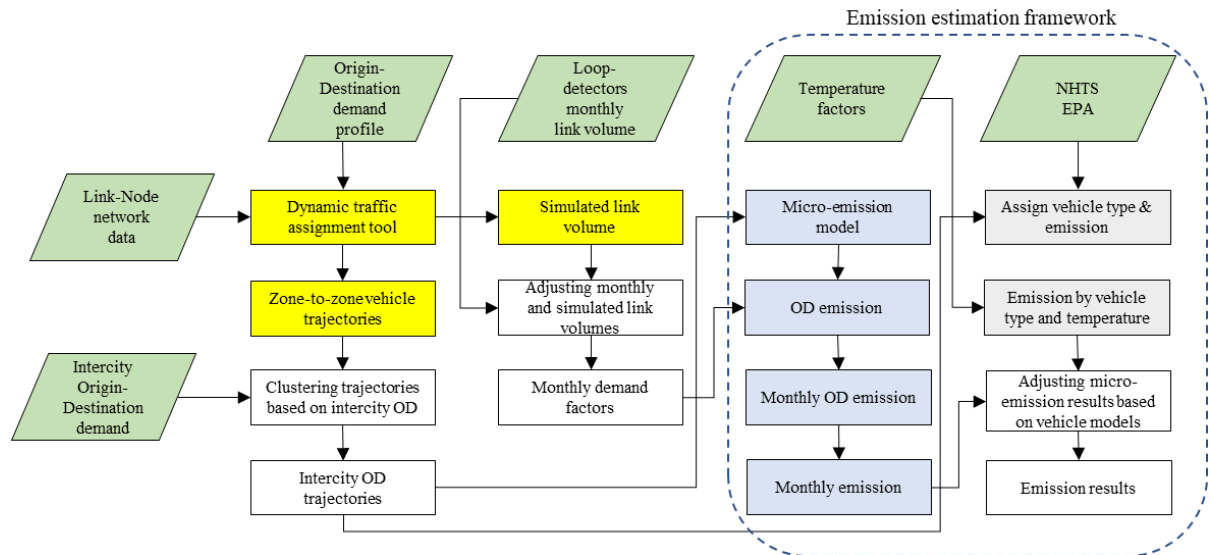

Figure 1 Different components of the proposed research framework to estimate emission in the intercity network of Michigan

### 3.1. Dynamic Traffic Simulation

A mesoscopic dynamic traffic simulation tool is used in this study to simulate the movement of individual vehicles throughout the state-wide network of Michigan. The static daily OD demand table for this network is provided by the Michigan Department of Transportation (MDOT). This table is converted to a time-dependent OD demand using time factors calculated based on observed hourly counts provided by installed loop detectors on Michigan highways to capture the traffic dynamic throughout the day [13]. This time-dependent demand, along with the road network properties, are the inputs to the dynamic traffic simulation tool used in this study, DYNASMART. This tool is used to simulate traveler's route choice behavior and find the network user equilibrium by assigning the vehicles to their shortest paths iteratively. This simulation tool provides the daily simulated vehicle trajectories, as well as information regarding the instantaneous speed, acceleration, and position in the network for the base scenario.

The zone-to-zone vehicle trajectories of the state-wide network include urban and intercity trajectories. As this study aims to study intercity trips and explore the effects of EV adoption on emission, the urban trajectories are excluded. To this end, the zones within each major city, which has a population of over 50,000, are clustered as an intercity node. The trajectories between the intercity nodes are considered as the intercity trajectories. The intercity network is discussed in detail in section 4.1.

### 3.2. Micro-Emission Model

The micro-emission model implemented in this study is a nonlinear regression model developed by Panis et al. [16]. This model considers the trip dynamics (speed and acceleration) and estimates different results for vehicles based on their category. It considers five different vehicle categories, including petrol car, diesel car, liquefied petroleum gas car, heavy-duty vehicles, and bus. The model is formulated as:

$$E_n^p(t) = max\left[0, c_1^p + c_2^p v_n(t) + c_3^p v_n^2(t) + c_4^p a_n(t) + c_5^p a_n^2(t) + c_6^p v_n(t)a_n(t)\right] \quad (1)$$

In which, $E_n^p(t)$ is the emission of pollutant $p$ (gram/second) from vehicle type $n$ at time $t$ (second). $v_n(t)$ and $a_n(t)$ are the speed $(m/s)$, and acceleration $(m/s^2)$, respectively, at that specific instant. The values $c_1^p$, …, $c_6^p$ are the nonlinear model coefficients. This model differentiates between acceleration and deceleration for certain pollutants, i.e., it has different coefficients for different pollutants. Since the scope of this study is the impact of light-duty vehicle fleet electrification on emission, calibrated models for the petrol car are adopted to estimate $CO_2$ emissions.

### 3.3. Monthly Travel Demand Factors

The static OD demand tables are mainly estimated by planning agencies for a typical day in the fall, which is assumed to be representative of the average daily travel pattern throughout the year. While this assumption is valid for urban trips, the travel demand usually fluctuates significantly throughout the year for intercity trips. Exploring the impacts of vehicle electrification on emission requires considering these fluctuations. In this study, these fluctuations are captured through OD-specific monthly travel demand factors estimated using the 24/7 count data from continuous counting stations installed on Michigan roads provided by MDOT (See Figure 2 for total observed monthly counts in all stations). Assigning the OD demand to the intercity network, the simulated link flows are iteratively calculated based on the OD factors to match the data from counting stations, assuming a proportional relationship between counting stations data and OD demands. From the 122 continuous counting stations provided by MDOT, 66 of them are used to capture the OD demand factors, as they are located on the intercity network links. For more information, please refer to Fakhrmoosavi et al. [15].

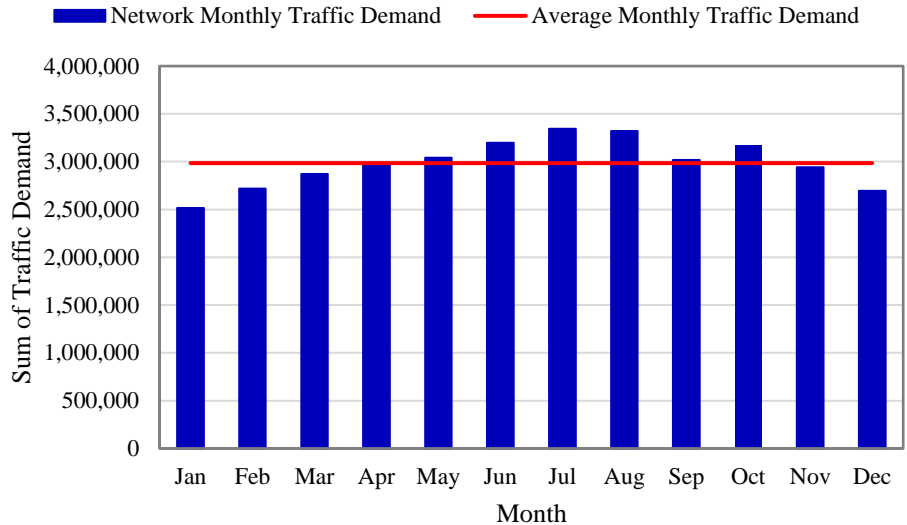

Figure 2 Total monthly demand for Michigan

### 3.4. Vehicle Type Adjustment Factor for Emission Estimation

As this study aims to estimate the light-duty vehicle emission and explore emission savings due to electrification, it is essential to distinguish between different vehicle types to estimate the emission accurately. Thus, information on the actual distribution of vehicle types defined by vehicle make, model, and manufactured year is required. This study uses the 2017 National Household Travel Survey (NHTS) sample dataset to estimate the distribution of various types of light-duty vehicles in the state of Michigan, to investigate the impact of vehicle type on the emission rate. This sample dataset includes vehicles categorized by vehicle make, vehicle model, and model year. This dataset is coupled with the associated $CO_2$ emission of each vehicle type, derived from the U.S. EPA dataset [40], as the main running emission of light-duty vehicles. The comparative emission rates per unit of distance for different vehicle types provide adjustment factors to be applied to the vehicles in different categories.

Comparing the NHTS dataset to the EPA dataset, only the vehicles that their $CO_2$ emission was available in the NHTS dataset were considered in this study. Based on this subset of NHTS data for Michigan, 570 different types of light-duty gasoline vehicles (almost 4,242,000 vehicles) are available from model years 2000 to 2017 [44]. Figure 3(a) shows the share of vehicles for each model year in Michigan based on NHTS data. The weighted average age of the vehicles in 2020 is 10.2 years, i.e., the average model year is 2011. Also, the weighted standard deviation of vehicles' age is 4.2 years. Figure 3(b) shows the share and distribution of vehicles in Michigan based on $CO_2$ emission rates. The weighted average of $CO_2$ emission rates over all vehicles in Michigan is 246 gr/mile based on NHTS. The weighted standard deviation of the $CO_2$ emission rate for vehicles in Michigan is 55 gr/mile.

EPA dataset provides vehicle emission rates based on make, model, engine type, and model year. The emission rates include the vehicle test results at the EPA lab in Ann Arbor, Michigan, as well as the data submitted by car companies [40]. The EPA dataset includes more detail than NHTS dataset because it entails results for different engines and different tests. Therefore, to make these two datasets comparable, for each vehicle make, model, and model year, the different emission rates of different tests and engine types are averaged in the EPA dataset to reach to the level of aggregation in NHTS dataset.

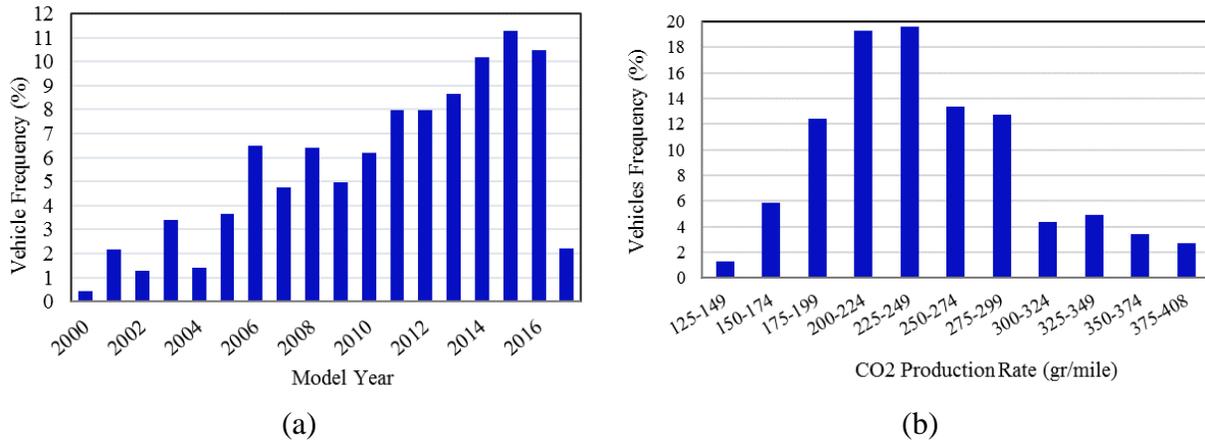

(a) (b)

Figure 3 Share of vehicles in Michigan for different vehicle specifications based on NHTS for (a) Model year (b) $CO_2$ production

To estimate the state-wide emission, vehicle types in the NHTS dataset are randomly assigned to the trajectories derived from the dynamic traffic assignment proportional to their market shares in Michigan. Then, the selected micro-emission model is employed to estimate the $CO_2$ emission for each trajectory. Since the implemented micro-emission model is calibrated based on a group of sedan vehicles with model years of 2005 or earlier [16], it does not accurately estimate other vehicle types' emissions depending on their body style and model year. Therefore, the micro-emission model estimations are adjusted to capture the impact of the vehicle type. To this end, the gasoline vehicles used in the calibration of the selected micro-emission model are selected and crossed with available vehicle emission data from the EPA dataset. The common vehicles between these two datasets are Toyota Corolla, Toyota Celica, and Volkswagen Golf. The average $CO_2$ emission rate of these vehicles is considered as the average emission rate of the selected micro-emission model. Thus, the emission factor of each vehicle type, $\mu_i$, is defined as below:

$$\mu_t = \frac{E_t}{E_b} \qquad (1)$$

Where $E_t$ is the EPA emission rate for vehicle type $t$ depending on its model year and make and $E_b$ is the average emission rate of the selected micro-emission model. Thus, the micro-emission model results based on the vehicle trajectories are adjusted by applying each vehicle type's emission factor.

### 3.5. Temperature Adjustment Factor for $CO_2$ Emission

While most of the emission estimation experiments are conducted at normal temperatures (70˚F), long and harsh winters, similar to what Michigan experiences frequently during the wintertime, necessitate the emission estimation in subzero temperatures. Therefore, this study investigates the impacts of temperature on emission rates to account for harsh weather on $CO_2$ emission rates by calculating another adjustment factor.

To capture the effects of temperature variations on $CO_2$ emission rates, Suarez-Bertoa and Astorga carried out a study on a group of gasoline vehicles with engine displacements of 1,000 to 2,000 ccs at the ambient temperatures of 20°F and 73°F [39]. The results of these tests are provided in Figure 4, where the dashed line shows the average $CO_2$ emission rate for the vehicles used to calibrate the Panis micro-emission model. Based on this study, on average, the $CO_2$ emission rate of vehicles at 20°F is 1.11 times more than the $CO_2$ emission rates at 73°F. This ratio is defined as the temperature adjustment factor in this study. This factor is applied to the emission calculated from the Panis micro-emission model for monthly emission estimation.

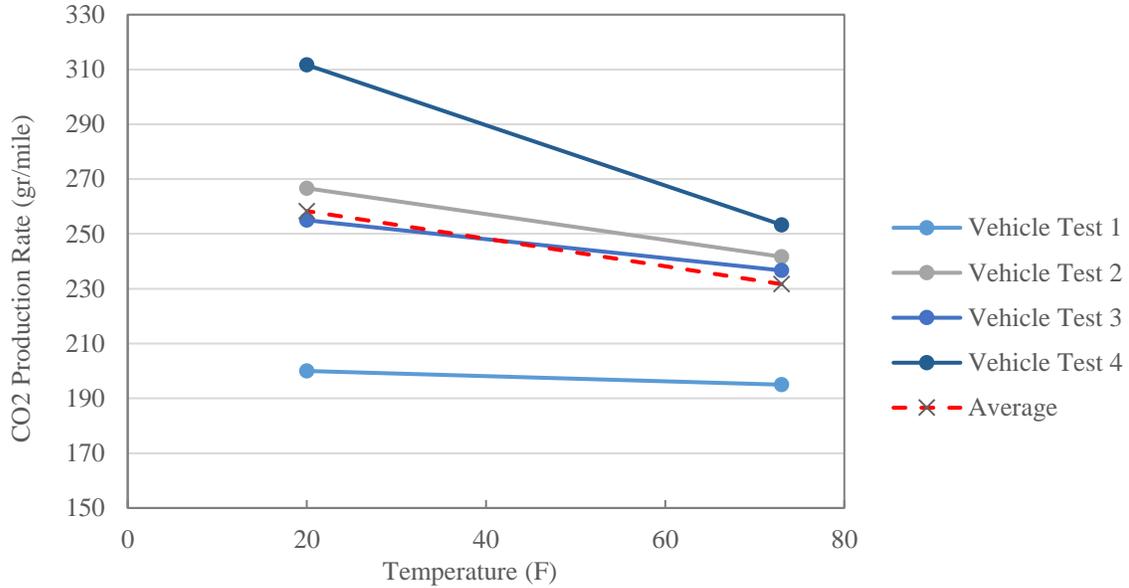

Figure 4 $CO_2$ production rates for gasoline vehicles under different ambient temperatures

### 3.6. Macro-Emission Estimation

The macro-emission estimation approach, used in this study as a benchmark, estimates emission only using the vehicle miles traveled in the network and the EPA average emission rate per unit of distance. The dynamic traffic assignment tool provides the trajectories of each trip, including the trip distances. The average vehicle $CO_2$ emission rate, $(\gamma)$, is estimated as 368 (gr/mile) by EPA [40]. Therefore, the daily macroscopic emissions, $e^i$, regardless of ambient temperature and vehicle types can be calculated for each OD pair $i$.

$$e^i = \gamma \sum_{n=1}^{N} D_n^i \qquad (2)$$

Where $N$ is the number of daily trajectories in the OD pair $i$, and $D_n^i$ is the trip distance of $n^{th}$ trip for that OD pair. Based on the number of days and the temperature factor for each month, the emissions can be adjusted as follows:

$$E_m^i = e^i N_m T_m \varphi_m \qquad (3)$$

Kavianipour et al.                                                                                                                                    10Where the $E_m^i$ is the adjusted monthly macro-emission of OD pair $i$, $N_m$ is the number of days in month $m$, $T_m$ is the temperature adjustment factor of that month, and $\varphi_m$ is the monthly demand factor. The results of the macro- and micro- emission models considering different adjustment factors (monthly demand factors, vehicle-type adjustment factor, and the temperature adjustment facto) are compared in the numerical result sections for the Michigan intercity network.

## 4. NUMERICAL EXPERIMENTS

In this section, the proposed methodology is implemented to estimate $CO_2$ savings resulting from the electrification of intercity trips in Michigan. First, the specifications of the traffic simulation model and the proposed intercity network are presented. Next, the vehicle emission estimation is carried out considering different adjustment factors and compared with the benchmark VMT-based approach. Then, a sensitivity analysis is conducted to demonstrate the impacts of different adoption scenarios in terms of vehicle types to be replaced by EVs. Next, the required infrastructure investment costs and the electricity demand for electrification of intercity trips are compared with the emission savings for different projected market shares of EVs. Finally, an economic analysis is conducted to explore the benefit to cost ratio of Michigan intercity network electrification.

### 4.1. Study Area

The network of interest is the road network of Michigan, provided by the MDOT, and shown in Figure 5. This network includes 2,330 traffic analysis zones (TAZs) and almost 30,000,000 daily trips as travel demand for a weekday in fall underlying normal weather conditions. To consider the intercity trips, the TAZs are clustered into 23 nodes representing major cities, which have populations over 50,000. This aggregation results in the intercity network that includes trajectories between these 23 nodes as intercity trajectories, representing almost 3,000,000 daily trips. The configuration of the intercity network is demonstrated in Figure 5.

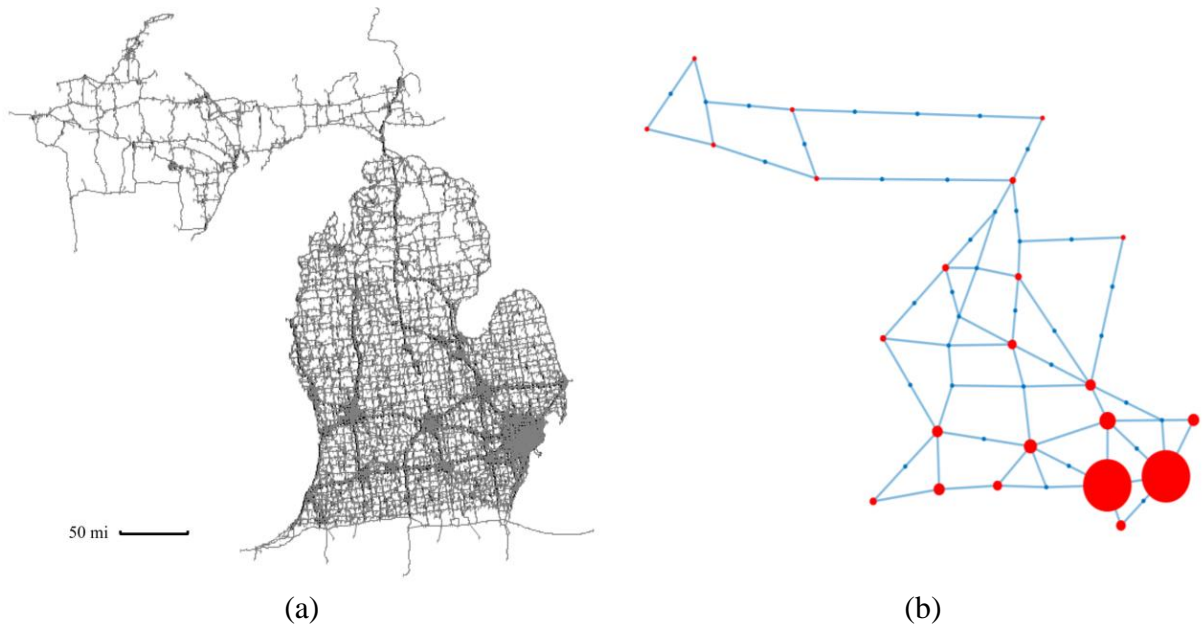

(a)                                                                                     (b)

Figure 5 (a) Transportation network of Michigan (b) Intercity network of Michigan

### 4.2. Vehicle Emissions

Micro-emission estimation models are calibrated to calculate emissions based on travel pattern and driving behavior. In this study, the micro-emission estimation model, presented by Panis et al. [16], is implemented to estimate emission for each daily intercity trajectory based on the speed and acceleration profiles, which are the outputs of the mesoscopic dynamic traffic assignment tool, DYNASMART. Based on the number of days in each month, temperature adjustment factor, and demand factor associated with each month, the adjusted monthly $CO_2$ emission can be calculated for any given OD pair with positive intercity demand. Then, the monthly $CO_2$ emission associated with intercity trips in Michigan can be calculated by summing up the estimated emission for all OD pairs, as presented in Table 1. In this table, the VMT associated with each month is calculated based on the vehicle miles traveled by all intercity trajectories and monthly demand factors. For each month, the $CO_2$ emission is calculated based on the EPA dataset and Panis micro-emission model, which are adjusted by the Michigan vehicle types and temperature factors. Applying the vehicle type adjustment factors reduces the estimated emission by the EPA general average method (EPA average $CO_2$ per unit of distance multiplied by total VMT), which demonstrates the average vehicle considered by the EPA generates more emission relative to vehicles on Michigan roads. Also, the Panis micro-emission model is calibrated based on efficient light sedan vehicles, which results in an underestimation. In addition, for both approaches (distance-based and micro-emission-based), considering the cold ambient temperature increases the estimated emission. Finally, the distance-based approach underestimates the generated emission relative to the micro-emission model due to ignoring the vehicle dynamics (speed and acceleration variations).

TABLE 1 Estimated $CO_2$ Emission Using a Variety of Approaches

| Month | VMT (million miles) | $CO_2$ Emission (million ton) | | | | | | | |
|---|---|---|---|---|---|---|---|---|---|
| | | EPA | | | | Micro-Emission Model | | | |
| | | General Average | Adjusted by Temperature | Michigan Vehicle Types | Michigan Vehicles-Adjusted by Temp | Base Value | Adjusted by Temperature | Michigan Vehicle Types | Michigan Vehicles-Adjusted by Temperature |
| Jan | 2,407 | 0.89 | 0.98 | 0.59 | 0.66 | 0.59 | 0.65 | 0.75 | 0.83 |
| Feb | 2,317 | 0.85 | 0.95 | 0.57 | 0.63 | 0.57 | 0.63 | 0.72 | 0.80 |
| Mar | 2,727 | 1.00 | 1.11 | 0.67 | 0.75 | 0.66 | 0.74 | 0.84 | 0.94 |
| Apr | 2,721 | 1.00 | 1.00 | 0.67 | 0.67 | 0.66 | 0.66 | 0.84 | 0.84 |
| May | 3,019 | 1.11 | 1.11 | 0.74 | 0.74 | 0.73 | 0.73 | 0.92 | 0.92 |
| Jun | 3,115 | 1.15 | 1.15 | 0.77 | 0.77 | 0.75 | 0.75 | 0.95 | 0.95 |
| Jul | 3,357 | 1.24 | 1.24 | 0.83 | 0.83 | 0.80 | 0.80 | 1.01 | 1.01 |
| Aug | 3,348 | 1.23 | 1.23 | 0.82 | 0.82 | 0.80 | 0.80 | 1.02 | 1.02 |
| Sep | 3,047 | 1.12 | 1.12 | 0.75 | 0.75 | 0.73 | 0.73 | 0.93 | 0.93 |
| Oct | 3,092 | 1.14 | 1.14 | 0.76 | 0.76 | 0.75 | 0.75 | 0.95 | 0.95 |
| Nov | 2,783 | 1.03 | 1.03 | 0.69 | 0.69 | 0.67 | 0.67 | 0.85 | 0.85 |
| Dec | 2,557 | 0.94 | 1.05 | 0.63 | 0.70 | 0.62 | 0.69 | 0.79 | 0.88 |
| Annual | 34,490 | 12.71 | 13.11 | 8.50 | 8.77 | 8.32 | 8.59 | 10.57 | 10.91 |

Figure 6 represents the monthly $CO_2$ emission estimated by three approaches. It demonstrates the fluctuations of the EPA general average method, the Panis micro-emission model, and the modified Panis micro-emission model considering the vehicle type and temperature adjustment factors. This figure shows that considering the Michigan vehicle type composition increases the base values of $CO_2$ emissions estimated by the Panis micro-emission model. This is mainly due to the fact that the Panis micro-emission model is calibrated using a group of sedan vehicles that vary in type and level of emission with the vehicle type composition in Michigan. Figure 6 also shows that EPA overestimates the $CO_2$ emission throughout the year. It is worth noting that even though $CO_2$ emission rates increase because of cold temperatures, the increased VMT during the summer months cancels it out, and the total $CO_2$ emission increases during the warmer months in Michigan. Note that the EPA average method overestimates the emission by 16 percent and the adjusted EPA average method underestimates the emission by 20 percent relative to the proposed methodology in this study, which includes incorporating a micro-emission model coupled with a mesoscopic traffic simulation tool while considering adjustment factors for vehicle type composition and temperature.

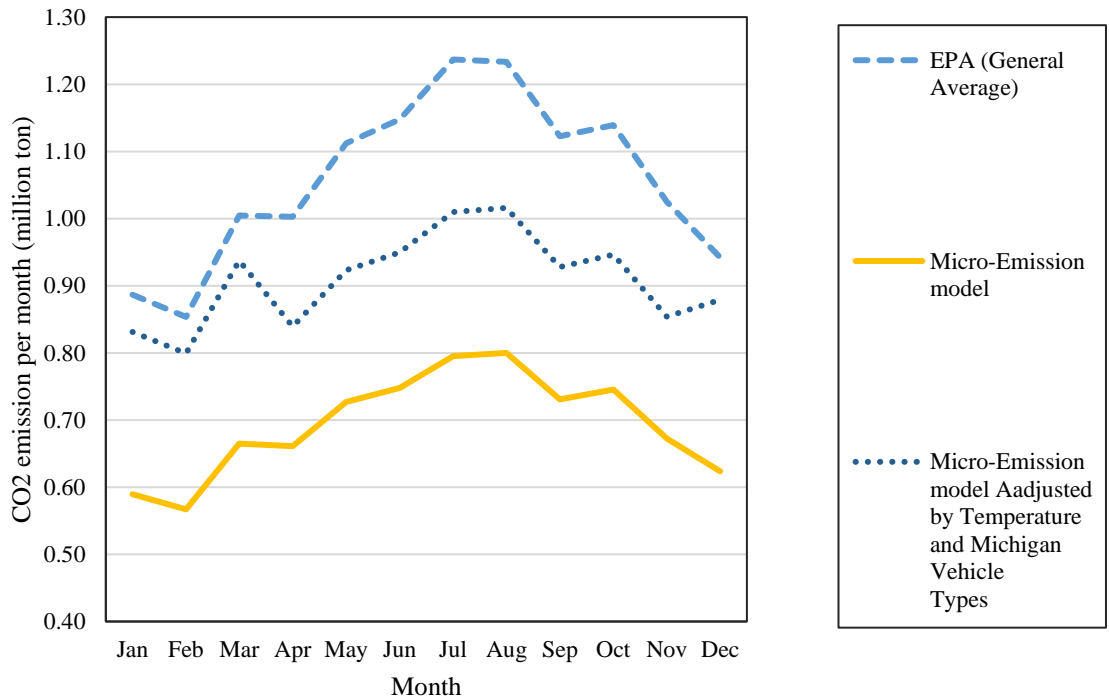

Figure 2 Monthly $CO_2$ emission calculated using different approaches

### 4.3. Scenario Analysis

Due to the uncertainty in the future of EV adoption, it is important to define possible scenarios that can capture the future vehicle choice of users and test the effect of each scenario on $CO_2$ emission reduction. Thus, in this section, different adoption scenarios are introduced for the vehicle types which will be replaced by EVs in future. The impact of these scenarios is investigated for projected Michigan EV market share in 2030. Then, the required charging infrastructure and energy demand for different EV market shares and technology scenarios are investigated. Finally, an economic analysis is conducted to show the profitability of Michigan intercity road network electrification considering the required investment to build the charging infrastructure.

#### 4.3.1. Adoption Scenarios

The type of vehicles that will be replaced by EVs affects the emission saving associated with the EV adoption. This choice can also reflect the governments' policies and subsidies to retire a certain type of vehicles and replace them with EVs. Therefore, five different scenarios are considered for the EV adoption. The first scenario, Random, considers a random selection of trajectories with a uniform distribution to be replaced by EVs on the road network. The second scenario, Old-Random, assumes that EVs will only replace vehicles older than 10-years old that are chosen randomly to match the projected market shares in future. This assumption considers the fact that even though in some households a vehicle newer than 10-years old is replaced by an EV, the sold gasoline vehicle will replace an older vehicle in the market eventually. The third scenario, Old-Pessimistic, assumes that EVs will replace vehicles over ten years old that have the least $CO_2$ emission among the vehicles of their age. In the fourth scenario, Old-Optimistic, it is assumed that EVs will replace vehicles older than ten years that have the worst $CO_2$ emission among vehicles of the same age. This assumption may promote public policies providing incentives to vehicles with the most $CO_2$ emission in order to meet specific emission saving goals. In the last scenario, Oldest, the oldest vehicles in the market are replaced by EVs. This assumption also reflects policies supported by governments toward decreasing the vehicles $CO_2$ emission.

Given the projected 6% market share of EVs in 2030 [45], the $CO_2$ emission reduction for each adoption scenario is presented in Figure 7. This Figure shows the $CO_2$ savings on the left axis and its associated societal cost on the right axis. Considering a social cost of $50/ton for $CO_2$ [46], there is a linear relation between $CO_2$ emission and its societal cost. According to Figure 7, the Old-Pessimistic scenario results in the least $CO_2$ reduction (0.45 to 0.56 million-ton), while the Old-Optimistic scenario decreases the $CO_2$ emission the most (0.75-0.93 million-ton (about 8.5 percent)). Overall, the electrification of the intercity trips in Michigan may have a range of savings of 0.56-0.93 million-ton $CO_2$ emission annually considering various replacement scenarios.

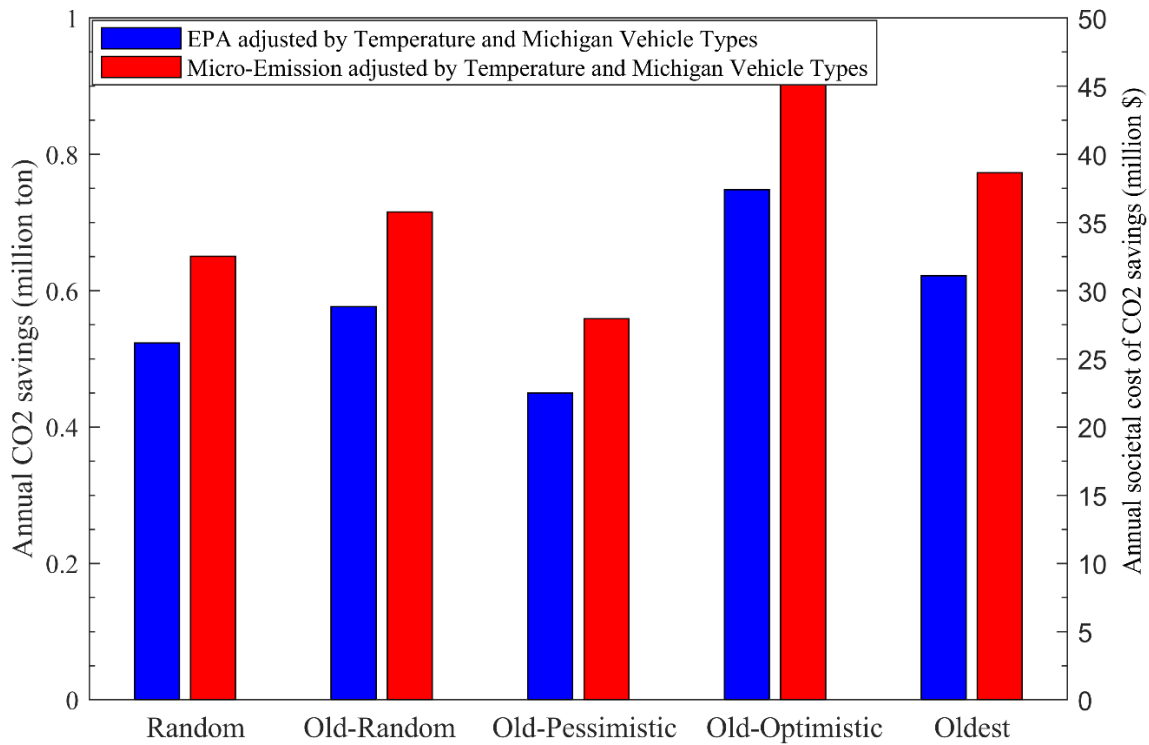

Figure 3 The $CO_2$ emission reduction for five scenarios with EV market share of 6 percent in 2030

### 4.3.2. **Market Share Analysis**

This section investigates the impacts of different EV market shares on the $CO_2$ savings in Michigan considering the source of electricity. To fully realize the impacts of the fleet electrification on the emission reduction, the $CO_2$ emission associated with the source of electricity needs to be considered. Although EVs remove the tailpipe emission, the electricity generation needs to be "Green" to effectively reduce the system emission. While the electricity mix consists of different sources, e.g., nuclear, wind, solar, and coal, the peak-hour use is generally supported by fossil fuel exhaustion. Providing a network of DC fast chargers associated with large-sized batteries can mitigate the network peak electricity demand and shift the energy demand to off-peak, e.g., nighttime.

The projected market shares of EVs in Michigan in 2025 and 2030 are 3 percent and 6 percent, respectively [45]. In addition to these market shares, a 10 percent EV market share is also explored in this section as an optimistic adoption of EVs. Assuming a battery performance of 3.5 (mile/kWh), the required charging infrastructures to support EV intercity trips in Michigan are calculated based on a recent study by Ghamami et al. [14]. For each market share, three technology scenarios are considered, namely Low-Tech, Mixed-Tech, and High-Tech. The Low-Tech scenario considers a battery capacity of 70 kWh for all EVs and 50 kW chargers for all chargers. The Mixed-Tech scenario considers a battery of 70 kWh and chargers with charging power of 150 kW. The High-Tech scenario considers a 100-kWh battery and 150 kW power for EVs and chargers, respectively. The technology scenarios are defined through a series of stakeholder meetings with charging station companies and car manufacturers [45]. The required charging infrastructure and energy demand at intercity chargers is provided in Table 2 for each technology and market share scenario for Michigan. The infrastructure cost represents the cost of charging stations, chargers, land acquisition, electricity provision, and maintenance over the lifetime of charging facilities, which is 10 years. (For more information on how to estimate the required charging infrastructure please refer to Ghamami et al., and Kavianipour et al. [14,42]). This table shows that deploying chargers with higher power can decrease the required charging infrastructure cost. Further, it shows that the required energy is more dependent on the battery size rather than the charging power. In the High-Tech scenarios, where the battery capacity is 100 kWh, the energy demand is 40 percent of the other two technology scenarios. Therefore, Low-Tech and Mixed-Tech scenarios are not only more expensive than the High-Tech scenario, but also put more energy demand on the intercity charging infrastructure.

Table 2 Required charging infrastructure and energy demand for different market shares and technology scenarios

| EV Market Share (%) | Technology Scenario | Number of Stations | Number of Chargers | Infrastructure Investment Cost for 10 years (million dollar) | Daily Charging Station Energy Demand (MWh) | Annual Charging Station Energy Demand (BWh) | Total Energy Demand (BWh) | Intercity Chargers Energy Demand Percent |
|---|---|---|---|---|---|---|---|---|
| 3 | Low-Tech | 34 | 250 | 14.30 | 69.73 | 25.45 | 295.62 | 8.61 |
| 3 | Mixed-Tech | 29 | 90 | 12.54 | 70.48 | 25.72 | 295.62 | 8.70 |
| 3 | High-Tech | 18 | 47 | 7.10 | 30.80 | 11.24 | 295.62 | 3.80 |
| 6 | Low-Tech | 38 | 478 | 23.07 | 140.04 | 51.12 | 591.24 | 8.65 |
| 6 | Mixed-Tech | 33 | 160 | 18.78 | 141.36 | 51.60 | 591.24 | 8.73 |
| 6 | High-Tech | 21 | 101 | 11.82 | 61.51 | 22.45 | 591.24 | 3.80 |
| 10 | Low-Tech | 41 | 748 | 33.10 | 235.22 | 85.85 | 985.40 | 8.71 |
| 10 | Mixed-Tech | 34 | 269 | 27.50 | 237.99 | 86.87 | 985.40 | 8.82 |
| 10 | High-Tech | 29 | 152 | 17.48 | 102.58 | 37.44 | 985.40 | 3.80 |

The total intercity miles traveled in the intercity network is 34,489 million miles according to the used travel model in this study provided by MDOT . Therefore, the total energy demand for each EV market share can be calculated based on the total miles traveled, market share, and battery performance, as provided in Table 2. Comparing the annual energy demand of EVs with the DC fast charging energy demand in Table 2, 3.8 to 8.8 percent of total energy needs to be provided from DC fast chargers

The $CO_2$ savings for each market share and adoption scenario is provided in Table 3. In this table, Adjusted EPA represents the values that are derived based on the EPA average data and are adjusted based on Michigan vehicle composition and temperature factors. Also, the Adjusted Micro-Emission shows the values that are calculated based on the Panis micro-emission model and are adjusted based on Michigan vehicle composition and temperature factors. This Table shows that $CO_2$ savings ranges from 0.21 million ton to 1.42 million ton annually depending on the EV market share and adoption scenario. The societal cost of $CO_2$ for each market share and adoption scenario can be calculated based on the linear relationship between the $CO_2$ savings and the social cost of $CO_2$ , which is $50/ton [46].

Table 3 $CO_2$ savings (million ton) for different adoption scenario and market shares in the intercity network of Michigan

| EV Market Share (%) | Random | | Old-Random | | Old-Pessimistic | | Old-Optimistic | | Oldest | |
|---|---|---|---|---|---|---|---|---|---|---|
| | Adjusted EPA | Adjusted Micro-Emission | Adjusted EPA | Adjusted Micro-Emission | Adjusted EPA | Adjusted Micro-Emission | Adjusted EPA | Adjusted Micro-Emission | Adjusted EPA | Adjusted Micro-Emission |
| 3 | 0.26 | 0.32 | 0.29 | 0.36 | 0.21 | 0.26 | 0.40 | 0.50 | 0.33 | 0.41 |
| 6 | 0.52 | 0.65 | 0.58 | 0.72 | 0.45 | 0.56 | 0.75 | 0.93 | 0.62 | 0.77 |
| 10 | 0.88 | 1.09 | 0.97 | 1.20 | 0.82 | 1.02 | 1.14 | 1.42 | 1.04 | 1.29 |

#### 4.3.3. Economic Analysis

This section compares the required charging infrastructure investments with the societal cost of $CO_2$ savings associated with different EV market shares for Michigan. To this end, the Old-Random scenario is considered as the future adoption scenario and the $CO_2$ savings of this scenario is compared with the charging infrastructure investment for EV market shares of 3, 6, and 10 percent for different technology scenarios, as presented in Table 4. The infrastructure cost in this table is converted to the depreciation cost per year considering a lifetime of 10 years and a zero-inflation rate. It should be noted that the $CO_2$ savings are directly related to the reduced gasoline consumption in vehicles because of their replacement with EVs and do not include the life-cycle $CO_2$ of gasoline and electricity. Table 4 shows that increasing the EV market share or improving the technology of batteries and chargers increases the benefit to cost ratio. The benefit to cost ratio shows that no matter which technology or market share is selected, the economical savings of the $CO_2$ reduction always justify the electrification of the Michigan intercity road network in terms of tailpipe emission.

Table 4 Economic analysis of intercity network electrification

| EV Market Share (%) | Technology | Annual Societal Cost of $CO_2$ Savings (million dollar) | Annual Charging Infrastructure Cost | Benefit to Cost Ratio |
|---|---|---|---|---|
| 3 | Low-Tech | 17.85 | 1.43 | 12.48 |
| 3 | Mixed-Tech | 17.85 | 1.25 | 14.23 |
| 3 | High-Tech | 17.85 | 0.71 | 25.14 |
| 6 | Low-Tech | 35.79 | 2.31 | 15.51 |
| 6 | Mixed-Tech | 35.79 | 1.88 | 19.05 |
| 6 | High-Tech | 35.79 | 1.18 | 30.28 |
| 10 | Low-Tech | 59.84 | 3.31 | 18.08 |
| 10 | Mixed-Tech | 59.84 | 2.75 | 21.76 |
| 10 | High-Tech | 59.84 | 1.75 | 34.24 |

## 5. SUMMARY AND CONCLUSION

EVs have zero tailpipe emission and their adoption decreases fossil fuel dependency. As EV adoption can be promoted by providing adequate charging infrastructure, a recent study has developed an optimization tool to support the EV intercity trips via a network of DC fast chargers [14]. This optimization tool is employed to find the locations and required number of chargers at each location to support the projected EV market share of Michigan, 6 percent, in 2030 [45]. In this study, we estimate the emission savings in the intercity network of Michigan, resulting from different market shares for EVs using the optimization tool developed by Ghamami et al. [45]. The emission savings are then compared with the EV electricity consumptions along intercity corridors, and more specifically, at the proposed charging stations located by the optimization model presented in Ghamami et al. [45]. Emission savings are a function of total intercity network emission and EV market share. To estimate the total intercity network emission, a state-of-the-art framework is developed. This framework first incorporates a micro-emission model to estimate the emission associated with each simulated vehicle trajectory. Then, it uses temperature and vehicle type adjustment factors to adjust the emission of each vehicle trajectory. Finally, the aggregated emission over vehicle trajectories is adjusted by monthly demand factors to calculate the annual network-wide emission. The estimated network-wide emission is then compared with a benchmark approach that incorporates annual VMT and average emission rate per unit of distance provided by EPA. Once the network emission is properly estimated, vehicle trajectories are randomly assigned to EVs to estimate emission savings. The random assignment of the vehicle trajectories to EVs is performed based on five different scenarios. The main findings of this study can be summarized as follows:

- A unique framework is presented and implemented for Michigan intercity trips that estimate emission at a large-scale network by coupling a micro-emission model with a mesoscopic traffic simulation tool considering monthly demand and temperature variations.

Kavianipour et al.                                                                 17

- The micro-emission model needs to be adjusted to account for the existing vehicle type composition and temperature variations in Michigan. Neglecting these factors results in a significant misestimation of emission savings.
- The EPA average method, incorporating VMT, overestimates the emission savings due to neglecting vehicle type composition specific to Michigan. Even the adjusted EPA average method, which considers vehicle type composition and temperature variations, underestimates the emission saving due to neglecting travel dynamics.
- The electrification of the intercity road network in Michigan can save, on average, 0.56-0.93 million tons of $CO_2$ emission annually in 2030. Thus, for a 6 percent EV market share, the annual intercity $CO_2$ emission will be reduced by 5.3 to 8.5 percent.
- In 2030, only 3.8 to 8.7 percent of the total EV energy requirements must be provided by DC fast chargers, which is between 22.45 to 51.60 BWh. The rest of the required energy to operate EVs on intercity corridors is supplied by home/ overnight charging.
- For a 3 percent EV market share, the annual $CO_2$ savings range from 0.26 to 0.5 million ton, while they reach to 1.02 to 1.42 million to for a 10 percent EV market share.
- The societal cost of $CO_2$ savings to charging infrastructure cost ratio varies between 12.5 to 34 depending on the technology scenario and EV market share.
- The societal cost of $CO_2$ savings ratio to charging infrastructure cost estimated in this study for different scenarios can be used by city planners and policy makers to justify the investment in the network electrification and users' incentive policies.

The presented framework in this study provides an important initial step towards exploring the emission savings resulting from fleet electrification. Future studies need to be conducted to explore the electricity source mix, which changes throughout the day with the variations in electric grid demand. In addition, this study focuses on intercity trips, and more research is needed to explore emission savings and required infrastructure investment for electrification of urban trips.

**ACKNOWLEDGMENTS**
This material is based upon work supported by the Department of Energy and Energy Services under Award Number EE008653. The authors also appreciate the assistance of the Bureau of Transportation Planning staff at the Michigan Department of Transportation (MDOT) in making data available to the study, especially Bradley Sharlow and Jesse Frankovich. The authors naturally remain solely responsible for all contents of the paper.

**AUTHOR CONTRIBUTIONS**
All authors contributed to all aspects of the study, from conception and design to analysis and interpretation of results, and manuscript preparation. All authors reviewed the results and approved the final version of the manuscript. The authors do not have any conflicts of interest to declare.

Kavianipour et al. 20